# LOCAL PARTIAL LIKELIHOOD ESTIMATION IN PROPORTIONAL HAZARDS REGRESSION

By Songnian Chen and Lingzhi Zhou

*Hong Kong University of Science and Technology*

Fan, Gijbels and King [*Ann. Statist.* **25** (1997) 1661–1690] considered the estimation of the risk function $\psi(x)$ in the proportional hazards model. Their proposed estimator is based on integrating the estimated derivative function obtained through a local version of the partial likelihood. They proved the large sample properties of the derivative function, but the large sample properties of the estimator for the risk function itself were not established. In this paper, we consider direct estimation of the relative risk function $\psi(x_2) - \psi(x_1)$ for any location normalization point $x_1$. The main novelty in our approach is that we select observations in shrinking neighborhoods of both $x_1$ and $x_2$ when constructing a local version of the partial likelihood, whereas Fan, Gijbels and King [*Ann. Statist.* **25** (1997) 1661–1690] only concentrated on a single neighborhood, resulting in the cancellation of the risk function in the local likelihood function. The asymptotic properties of our estimator are rigorously established and the variance of the estimator is easily estimated. The idea behind our approach is extended to estimate the differences between groups. A simulation study is carried out.

**1. Introduction.** The Cox proportional hazards model is by far the most popular model in survival analysis. In the Cox model, the conditional hazard rate of a survival time, $T$, given the regressor vector, $X = x$, is modeled as

$$(1.1) \qquad \lambda(t|x) = \lambda_0(t) \exp\{\psi(x)\},$$

where $\lambda_0(t)$ is the baseline hazard function and $\psi(x)$ is the risk function which measures the contribution of $X$ at $x$. Typically, the baseline hazard function is left unspecified, while the risk function $\psi(\cdot)$ is specified parametrically as $\psi(X) = X^T \beta^0$ for a vector of coefficients $\beta^0$, where $^T$ denotes the transpose of the vector. In practice, survival data are often censored due to termination of the study or early withdrawal from the study. We can thus observe only $Y = \min\{T, C\}$, where $C$ is the censoring time independent of $T$ given $X$. In addition, we also observe the censoring indicator, $\delta = I\{T \leq C\}$, as well as the covariate vector $X$. Let $\{(X_i, Y_i, \delta_i), i = 1, \ldots, n\}$ represent an i.i.d. sample from the population $(X, \min(T, C), I\{T \leq C\})$, $t_1^0 < \cdots < t_N^0$ denote the ordered observed failure times and $(j)$ provide the label for the item failing at $t_j^0$. Define $R_j$ as the risk set at time $t_j^0$: $R_j = \{i : Y_i \geq t_j^0\}$. Cox [5] suggested that estimation and inference on $\beta^0$ be based on the partial likelihood function

$$L(\beta) = \prod_{j=1}^{N} \frac{\exp(X_{(j)}^T \beta)}{\sum_{i \in R_j} \exp(X_i^T \beta)}.$$

With flexible specification for the baseline hazard function, Cox's model and the partial likelihood approach provide a very convenient way to measure the covariate effects on the survival time. See [7, 16, 19] and [21] for references on this model. However, the parametric specification of the risk factor $\psi(x)$ is assumed largely for convenience. In general, misspecification of the risk function will lead to inconsistent estimation and misleading statistical inferences. Therefore, it is desirable to relax the parametric specification of the risk function. When the risk function $\psi(x)$ is not parametrically specified, the partial likelihood function is of the form

$$L(\psi) = \prod_{j=1}^{N} \frac{\exp(\psi(X_{(j)}))}{\sum_{i \in R_j} \exp(\psi(X_i))}.$$

Note that the function $\psi(x)$ is only estimable up to a location normalization. In fact, in Cox's original setup with the linear specification for the risk function, no intercept is allowed in $\beta$. The relative risk provides all the information regarding the contribution of the covariates based on the Cox proportional hazards model. Recently, Fan, Gijbels and King [12] and Tibshirani and Hastie [23] considered local versions of the partial likelihood approach. They compared the relative risks for observations whose corresponding covariate values belong to a single shrinking neighborhood. Obviously, these observations all have the same risk factor to the first order. Namely, the first term of the Taylor expansion of $\psi(X_i)$ will be the same for $X_i$ close to $x$. As a result, their estimation is based on the second-order comparison of the relative risk factors, leading to the estimators for the derivatives of the function $\psi(x)$ only. To recover the original risk function, they suggested using the expression $\psi(x) = \int_0^x \psi'(t)\,dt$ by replacing $\psi'(t)$ with their local partial likelihood estimates $\hat{\psi}'(t)$ for $t \in (0, x)$. Note that $\psi(0) = 0$ is implicitly imposed for identifiability of $\psi(\cdot)$. However, the large sample property of the estimator $\int_0^x \hat{\psi}'(t)\,dt$ has not been formally established; thus, formal statistical inference is not feasible. Furthermore, their derivative-based approach cannot be extended to the case in which the covariate variable is discrete, since $\psi(\cdot)$ is canceled out in the local partial likelihood, as pointed out in Section 2.

In this paper, we consider the direct estimation of the relative risk function $\psi(x_2) - \psi(x_1)$ through a new version of the local partial likelihood. Intuitively, in constructing our local partial likelihood, we use observations in the neighborhoods of either $x_1$ or $x_2$, which have risk factors different to the first order, thus enabling direct estimation and inference of the relative risk $\psi(x_2) - \psi(x_1)$. Moreover, when the covariate variable is discrete, our approach reduces to the partial likelihood estimator for a two-sample comparison of survival times in the form of the usual proportional hazards model. In other words, our procedure reduces to an efficient estimation in the case of discrete covariates. Thus, we can expect our procedure to have high efficiency, at least when the data for the covariate variable $X$ are not evenly distributed.

Our procedure can be easily adapted to estimate differences in risk functions at any point $x$ between two different groups by constructing our local partial likelihood using observations in the neighborhood of $x$ for the two groups under consideration. We apply our procedure to the PBC data and find no treatment difference, which is consistent with findings from parametric hazard regression. In an attempt to compare treatment differences, Fan and Gijbels [11] suggested estimating $\psi'(x)$ for the two treatments separately and then integrating the derivative functions to recover the two risk functions. However, in this way each risk function is only estimable to within a constant. As a consequence, the treatment difference is not directly estimable. Fan and Gijbels [11] thus imposed zero risk for both treatments at the left endpoint of the support of $x$. Our numerical analysis found that this assumption is inappropriate.

There are some related studies on nonparametric regression techniques with censored data. Gentleman and Crowley [14] proposed an iterative algorithm to estimate $\psi(\cdot)$ with a uniform kernel function. Li and Doss [20]) investigated the nonparametric estimation of the conditional hazard and distribution functions using local linear fits. O'Sullivan [22], Hastie and Tibshirani [15] and Kooperberg, Stone and Truong [17, 18] used spline methods to study the model.

This paper is organized as follows. Section 2 introduces our estimator and the idea is extended to the case when discrete covariates are also present. A numerical study is presented in Section 3 where we compare our procedure with that of [12]. We also apply our procedure to the PBC data. Section 4 concludes the paper.

**2. Local partial likelihood estimators.** Recall that the partial likelihood function for model (1.1) is

$$(2.1) \qquad L(\psi) = \prod_{j=1}^{N} \frac{\exp(\psi(X_{(j)}))}{\sum_{i \in R_j} \exp(\psi(X_i))}$$

for an i.i.d. sample. For notational simplicity, we assume that $X$ is a continuously distributed random variable. We discuss cases with discrete and multivariate covariates later.

Suppose now that the form of $\psi(x)$ is not specified and that the $p$th-order derivative exists at $x$. Then a local model [11] of $\psi(X)$ can be expressed as

$$(2.2) \qquad \psi(X) \approx \psi(x) + \psi'(x)(X-x) + \cdots + \frac{\psi^{(p)}(x)}{p!}(X-x)^p$$

by Taylor expansion for $X$ in a neighborhood of $x$. Namely, $\psi(X) \approx \widetilde{X}^T \beta$ for $X$ close to $x$, where $\widetilde{X} = \{1, X-x, \ldots, (X-x)^p\}$ and $\beta = (\beta_0, \ldots, \beta_p)^T = \{\psi(x), \psi'(x), \ldots, \frac{\psi^{(p)}(x)}{p!}\}^T$. Let $K$ be a kernel function that smoothly downweighs the contribution of remote data points, let $h$ be the bandwidth parameter that controls the size of the local neighborhood and let $\widetilde{X}_i = \{1, X_i - x, \ldots, (X_i - x)^p\}^T$ for $i = 1, \ldots, n$. Fan, Gijbels and King [12] considered nonparametric estimation based on a local partial likelihood function

$$(2.3) \quad \sum_{j=1}^{N} K_h(X_{(j)} - x) \left[ \psi(X_{(j)}) - \log\left\{ \sum_{i \in R_j} \exp(\psi(X_i)) K_h(X_i - x) \right\} \right].$$

Using the local model (2.2), Fan, Gijbels and King [12] proposed to estimate the $\beta^*$ defined below with the likelihood function

$$(2.4) \quad \begin{aligned} &\sum_{j=1}^{N} K_h(X_{(j)} - x) \left[ \widetilde{X}_{(j)}^T \beta - \log\left\{ \sum_{i \in R_j} \exp(\widetilde{X}_i^T \beta) K_h(X_i - x) \right\} \right] \\ &= \sum_{j=1}^{N} K_h(X_{(j)} - x) \left[ \widetilde{X}_{(j)}^{*T} \beta^* - \log\left\{ \sum_{i \in R_j} \exp(\widetilde{X}_i^{*T} \beta^*) K_h(X_i - x) \right\} \right], \end{aligned}$$

where $K_h(t) = K(t/h)/h$,

$$\beta^* = (\beta_1, \ldots, \beta_p)^T \quad \text{and} \quad \widetilde{X}_i^* = \{X_i - x, \ldots, (X_i - x)^p\}^T.$$

Note, however, that the function value $\psi(x)$ is not directly estimable, as (2.4) does not involve the intercept $\beta_0 = \psi(x)$, which has been canceled out. Tibshirani and Hastie [23] considered a similar approach using a nearest neighborhood method. Furthermore, if $X$ is a discrete random variable taking on a finite number of values, a window around a value $x$ only contains that value itself. Therefore, (2.3) reduces to

$$\sum_{j=1}^{N} I(X_{(j)} = x)\left[\psi(x) - \log\left\{\sum_{i \in R_j} \exp(\psi(x))I(X_i = x)\right\}\right]$$

$$= \sum_{j=1}^{N} I(X_{(j)} = x)\left[-\log\left\{\sum_{i \in R_j} I(X_i = x)\right\}\right],$$

which no longer depends on $\psi(\cdot)$. This approach is thus not applicable to the case of discrete covariates.

2.1. *Estimation of the relative risk function.* We consider direct estimation of $\psi(x_2) - \psi(x_1)$ for a normalization point $x_1$ and any other point $x_2$ in the domain of $x$. By including observations in the neighborhoods of either $x_1$ or $x_2$, we consider the following local partial likelihood which essentially replaces $K_h(X_i - x)$ in (2.3) with $[K_h(X_i - x_1) + K_h(X_i - x_2)]$:

$$L_n = \sum_{j=1}^{N} [K_h(X_{(j)} - x_1) + K_h(X_{(j)} - x_2)]$$

(2.5)

$$\times \left[\psi(X_{(j)}) - \log\left\{\sum_{i \in R_j} \exp(\psi(X_i))[K_h(X_i - x_1) + K_h(X_i - x_2)]\right\}\right].$$

Let $\alpha = \psi(x_2) - \psi(x_1)$ and, for $l = 1, 2$,

$$\widetilde{X}_{li}^* = \{X_i - x_l, \ldots, (X_i - x_l)^p\}^T \quad \text{for } i = 1, 2, \ldots, n,$$

$$\beta_{x_l} = (\beta_{x_l,0}, \beta_{x_l,1}, \ldots, \beta_{x_l,p})^T$$

$$= \left\{\psi(x_l), \psi'(x_l), \ldots, \frac{\psi^{(p)}(x_l)}{p!}\right\}^T = (\beta_{x_l,0}, \beta_{x_l}^{*T})^T.$$

Using the local models in the neighborhoods of $x_1$ and $x_2$, we obtain an approximation of $L_n$,

(2.6) $$\tilde{L}_n = \tilde{L}_{n1} + \tilde{L}_{n2},$$

where

$$\tilde{L}_{n1} = \sum_{j=1}^{N} K_h(X_{(j)} - x_1)\widetilde{X}_{1(j)}^T \beta_{x_1}$$

$$- \sum_{j=1}^{N} K_h(X_{(j)} - x_1)\log\left\{\sum_{i \in R_j} [\exp(\widetilde{X}_{1i}^T \beta_{x_1})K_h(X_i - x_1)\right.$$

$$+ \exp(\widetilde{X}_{2i}^T \beta_{x_2}) K_h(X_i - x_2)]\bigg\}$$

(2.7)
$$= \sum_{j=1}^{N} K_h(X_{(j)} - x_1) \widetilde{X}_{1(j)}^{*T} \beta_{x_1}^*$$
$$- \sum_{j=1}^{N} K_h(X_{(j)} - x_1) \log\bigg\{\sum_{i \in R_j} [\exp(\widetilde{X}_{1i}^{*T} \beta_{x_1}^*) K_h(X_i - x_1)$$
$$+ \exp(\alpha + \widetilde{X}_{2i}^{*T} \beta_{x_2}^*) K_h(X_i - x_2)]\bigg\}$$

and, similarly,

$$\tilde{L}_{n2} = \sum_{j=1}^{N} K_h(X_{(j)} - x_2)(\alpha + \widetilde{X}_{2(j)}^{*T} \beta_{x_2}^*)$$

(2.8)
$$- \sum_{j=1}^{N} K_h(X_{(j)} - x_2) \log\bigg\{\sum_{i \in R_j} [\exp(\widetilde{X}_{1i}^{*T} \beta_{x_1}^*) K_h(X_i - x_1)$$
$$+ \exp(\alpha + \widetilde{X}_{2i}^{*T} \beta_{x_2}^*) K_h(X_i - x_2)]\bigg\}.$$

Clearly, our formulation of the local partial likelihood will allow direct estimation of $\alpha = \psi(x_2) - \psi(x_1)$, in contrast to [12] and [23].

In principle, one can estimate $\psi(x_2) - \psi(x_1)$ by finding the value that maximizes (2.6). Intuitively, however, observations in the neighborhood of $x_1$ would not be informative about the derivatives of $\psi(\cdot)$ at $x_2$, or vice versa. Therefore, a one-step estimator, which simultaneously estimates $(\alpha, \beta_{x_1}^*, \beta_{x_2}^*)$ through maximization of (2.6), is not particularly appealing. Instead, we adopt the following two-step strategy. In the first step, we adopt the approach of [12] to obtain the estimates $(\hat{\beta}_{x_1}^*, \hat{\beta}_{x_2}^*)$ for $(\beta_{x_1}^*, \beta_{x_2}^*)$. In the second step, we propose to estimate $\alpha$ by $\hat{\alpha}$, which maximizes

(2.9) $\quad \tilde{L}_n(\alpha, \hat{\beta}_{x_1}^*, \hat{\beta}_{x_2}^*) = \tilde{L}_{n1}(\alpha, \hat{\beta}_{x_1}^*, \hat{\beta}_{x_2}^*) + \tilde{L}_{n2}(\alpha, \hat{\beta}_{x_1}^*, \hat{\beta}_{x_2}^*).$

It is also worth noting that this approach is computationally more attractive. In addition, we find from our simulation that the performance of the two-step estimator is more stable than that of the one-step estimator.

Note that, from the construction of our estimator, only observations in the neighborhoods of either $x_1$ or $x_2$ will affect the estimation of $\psi(x_2) - \psi(x_1)$. The approach of [12], on the other hand, is cumulative in nature, in the sense that estimating $\psi(x_2) - \psi(x_1)$ requires the estimation of $\psi'(x)$ for $x$ in the interval between $x_1$ and $x_2$; consequently, a likely drawback of their approach is that inaccurate estimation of $\psi'(x)$ for $x$ in a neighborhood belonging to $[x_1, x_2]$ will adversely affect the precision in estimating $\psi(x_2) - \psi(x_1)$, which utilizes the estimates of $\psi'(x)$ for all of the $x \in [x_1, x_2]$. Indeed, this observation is confirmed in our simulation study presented in Section 3.1.

We now consider the asymptotic property of our estimator. Set

$$S(v|x) = P(Y \geq v|x).$$

We impose the following conditions:

1. the kernel function is a bounded symmetric density function with compact support;
2. the function $\psi(\cdot)$ has a continuous $(p_1 + 1)$st derivative around $x_1$ and $x_2$;
3. the density $f(\cdot)$ of $X$ is continuous at $x_1$ and $x_2$;
4. the conditional probability $S(v|\cdot)$ is equicontinuous at $x_1$ and $x_2$;
5. the local bandwidths $h$ and $h_1$ satisfy $h/h_1 \to 0$, $nh \to \infty$; $nh^{2p+3}$ and $nh_1^{2p_1+3}$ are both bounded, where $p_1 \geq p$ and $h_1$ are, respectively, the degree of polynomial and the bandwidth used for estimating the derivatives of [12].

The choice of bandwidths deserves some attention here. It is clear from (A.8) in the Appendix that when $h/h_1$ is not bounded, the asymptotic normality of $\hat{\alpha}$ cannot be achieved. When $h/h_1$ is bounded but does not approach 0 as $n \to \infty$, the asymptotic distribution of the first-step estimator will affect the distribution of $\hat{\alpha}$. When $h/h_1 \to 0$, which we impose here, the asymptotic distribution of the first-step estimator will have no impact on the second-step estimator, which makes the expression of the asymptotic variance much simpler than without the condition.

Our main result is stated in the following theorem.

THEOREM 1. *Under conditions 1–5, we have*

(2.10) $$\sqrt{nh}(\hat{\alpha} - \alpha - b_n(x_1, x_2)) \xrightarrow{D} N(0, \sigma^2(x_1, x_2)),$$

*where*

$$b_n(x_1, x_2) = \frac{h^{p+1}}{(p+1)!}\{\psi^{(p+1)}(x_2) - \psi^{(p+1)}(x_1)\}\left(\int u^{p+1} K(u)\, du\right) \quad \text{and}$$

$$\sigma^2(x_1, x_2) = \left(\int K^2(u)\, du\right)\left(\int_0^\infty \frac{a_{x_1}(v) a_{x_2}(v)}{a_{x_1}(v) + a_{x_2}(v)} \lambda_0(v)\, dv\right)^{-1},$$

*with $a_x(v) = e^{\psi(x)} f(x) P(v|x)$.*

The only unknown term of the bias, $\psi^{(p+1)}(x)$, can be easily estimated from our first-step estimator of the derivatives if we choose $p_1 > p$. To estimate $\sigma^2$, we observe that if we know $\psi$, then

$$\hat{a}_{x_l}(v) = n^{-1} \sum_{i=1}^{n} K_h(X_i - x_l) Y_i(v) \exp(\psi(X_i))$$

will have converged to $a_{x_l}(v)$ for $l = 1, 2$ and all $v$. In addition, the baseline hazard function $\Lambda_0(t)$ can be estimated by the Breslow estimator [3, 4],

$$\hat{\Lambda}_0(t) = \sum_{j=1}^{n} \frac{\delta_j I(Y_j \leq t)}{\sum_{i \in R_j} \exp(\psi(X_i))}.$$

Note that since $\psi(\cdot)$ itself is not estimable, $a_x(\cdot)$ or $\Lambda_0(\cdot)$ is only estimable up to scale. However, we can express $\int_0^\infty \frac{\hat{a}_{x_1}(v)\hat{a}_{x_2}(v)}{\hat{a}_{x_1}(v)+\hat{a}_{x_2}(v)} d\hat{\Lambda}_0(v)$ as

$$n^{-1} \sum_{j=1}^{n} \delta_j \frac{(\sum_{i \in R_j} K_h(X_i - x_1) \exp(D_i))(\sum_{i \in R_j} K_h(X_i - x_2) \exp(D_i))}{(\sum_{i \in R_j} \exp(D_i))(\sum_{i \in R_j} (K_h(X_i - x_1) + K_h(X_i - x_2)) \exp(D_i))},$$

where $D_i = \psi(X_i) - \psi(x_1)$ is already estimated with our methodology. Obviously, we can estimate $\sigma^2(x_1, x_2)$ by $\hat{\sigma}^2(x_1, x_2)$, where

$$\hat{\sigma}^2(x_1, x_2)$$
$$= \left(\int K^2(u) \, du\right)$$
$$\times \left\{ n^{-1} \sum_{j=1}^{n} \delta_j \left(\sum_{i \in R_j} K_h(X_i - x_1) \exp(\hat{D}_i)\right) \right.$$
$$\times \left(\sum_{i \in R_j} K_h(X_i - x_2) \exp(\hat{D}_i)\right)$$
$$\times \left(\left(\sum_{i \in R_j} \exp(\hat{D}_i)\right)\right.$$
$$\left.\left.\times \left(\sum_{i \in R_j} (K_h(X_i - x_1) + K_h(X_i - x_2)) \exp(\hat{D}_i)\right)\right)^{-1} \right\}^{-1},$$

with $\hat{D}_i$ an estimate of $D_i$.

The theoretical optimal bandwidth can be obtained by minimizing the asymptotic weighted mean integrated squared error

$$\int \int [\{b_n(x_1, x_2)\}^2 + \sigma^2(x_1, x_2)] w(x_1) w(x_2) \, dx_1 \, dx_2,$$

resulting in

$$h_{\text{opt}} = C_{0,p}(K) \bigg[ \int \int \bigg( \int_0^\infty \frac{a_{x_1}(v) a_{x_2}(v)}{a_{x_1}(v) + a_{x_2}(v)} \lambda_0(v)\, dv \bigg)^{-1}$$
$$\times w(x_1) w(x_2)\, dx_1\, dx_2$$
$$\times \bigg( \int \int (\psi^{(p+1)}(x_2) - \psi^{(p+1)}(x_1))^2$$
$$\times w(x_1) w(x_2)\, dx_1\, dx_2 \bigg)^{-1} \bigg]^{1/(2p+3)} n^{-1/(2p+3)},$$

where $C_{0,p}(K)$ are constants depending on $p$ and $K$. The value of $C_{0,p}(K)$ is tabulated in Table 3.2 of [11]. For a detailed discussion of the issue of model complexity, see [11].

REMARK 1. Note that when $X$ is a discrete random variable, our local partial likelihood (2.5) reduces to

$$L_n = \sum_{j=1}^N \bigg[ J_{1(j)} \psi(x_1) + J_{2(j)} \psi(x_2)$$
$$- (J_{1(j)} + J_{2(j)}) \log \bigg\{ \sum_{i \in R_j} (J_{1i} \exp \psi(x_1) + J_{2i} \exp(\psi(x_2))) \bigg\} \bigg]$$
$$= \sum_{j=1}^N \bigg[ J_{2(j)} \alpha - (J_{1(j)} + J_{2(j)}) \log \bigg\{ \sum_{i \in R_j} (J_{1i} + J_{2i} \exp(\alpha)) \bigg\} \bigg],$$

where $J_{li} = I(X_i = x_l), J_{l(j)} = I(X_{(j)} = x_l)$ for $l = 1, 2$, and $\alpha = \psi(x_2) - \psi(x_1)$. This is, in fact, the partial likelihood estimator for a two-sample comparison of survival time in the form of the proportional hazards model. In other words, our procedure yields an efficient estimation for the case of discrete covariates since it is well known that the partial likelihood estimator is efficient ([8] and [2]).

REMARK 2. Similar to the partial likelihood approach ([3] and [6]) and the local partial likelihood approach of [12] and [23], our version of the local partial likelihood function can also be viewed as a local profile likelihood. Analogous to [12], the local likelihood in our setting can be written as

$$\log L = \sum_{i=1}^n [\delta_i \{\log \lambda_0(Z_i) + \psi(X_i)\} - \Lambda_0(Z_i) \exp(\psi(X_i))]$$
(2.11)
$$\times (K_h(X_i - x_1) + K_h(X_i - x_2)).$$

Consider nonparametric modeling for $\Lambda_0(\cdot)$, which has a jump of $\lambda_j$ at $t_j$, $\Lambda_0(t, \lambda) = \sum_{j=1}^{N} \lambda_j I\{t_j \leq t\}$. Then

$$\Lambda_0(Z_i, \lambda) = \sum_{j=1}^{N} \lambda_j I\{i \in R_j\}.$$

Substituting these two expressions into the local likelihood expression (2.11), we obtain

$$\log L = \sum_{j=1}^{N} (K_h(X_i - x_1) + K_h(X_i - x_2))[\log \lambda_j + \psi(X_{(j)})]$$

(2.12)

$$- \sum_{i=1}^{n} \sum_{j=1}^{N} (K_h(X_i - x_1) + K_h(X_i - x_2)) \lambda_j I\{i \in R_j\} \exp(\psi(X_i)).$$

Then, by maximizing $\log L$ with respect to $\lambda_j$ $(j = 1, \ldots, N)$, we have

$$\hat{\lambda}_j = \frac{K_h(X_i - x_1) + K_h(X_i - x_2)}{\sum_{i \in R_j} (K_h(X_i - x_1) + K_h(X_i - x_2)) \exp(\psi(X_i))}.$$

Substituting $\hat{\lambda}_j$ into (2.12) yields

$$\max_{\lambda_0} \log L$$

$$= \sum_{j=1}^{N} (K_h(X_i - x_1) + K_h(X_i - x_2))$$

(2.13)

$$\times \left\{ \psi(X_{(j)}) - \log \sum_{i \in R_j} (K_h(X_i - x_1) + K_h(X_i - x_2)) \exp(\psi(X_i)) \right\}$$

$$+ \sum_{j=1}^{N} (K_h(X_i - x_1) + K_h(X_i - x_2))$$

$$\times \{ \log[K_h(X_i - x_1) + K_h(X_i - x_2)] - 1 \}.$$

Clearly, maximizing (2.13) is equivalent to maximizing (2.5).

REMARK 3. One referee pointed out that when estimating $\psi(x_2) - \psi(x_1)$, an alternative approach would be to estimate $\psi(x_2) - \psi(x_3)$ and $\psi(x_3) - \psi(x_1)$ separately and then to combine them for any point $x_3$ between $x_1$ and $x_2$. In general, these two approaches will lead to different estimates. Theoretical justification could be based on two separate asymptotic linear representations, as in (A.8). In our simulation experiment (not reported here), these two approaches seem to be comparable.

REMARK 4. One of the main advantages of the proportional hazards model is that it can easily accommodate time-varying covariates. The time-varying covariates can be incorporated into our approach in a straightforward way by expressing the local partial likelihood through the counting process representation. As our paper is largely based on its comparison with [12], we chose to use notation similar to that used in [12], in order to facilitate comparison.

2.2. *Estimating the differences between groups.* Our approach can be easily modified to accommodate more realistic situations in which estimating the differences between groups is necessary. Let the risk function be $\psi(x,z)$, where $x$ is continuous and $z$ is discrete, taking two values. We focus on estimating $\psi(x,z_2) - \psi(x,z_1)$. Similar to (2.5), we include observations in the neighborhoods of $(x, z_1)$ and $(x, z_2)$. Our local version of the log-likelihood becomes

$$L_n = \sum_{j=1}^{N} K_h(X_{(j)} - x)(I(Z_{(j)} = z_1) + I(Z_{(j)} = z_2))$$

$$\times \left[ \psi(X_{(j)}, Z_{(j)}) - \log \left\{ \sum_{i \in R_j} \exp(\psi(X_i, Z_i)) K_h(X_i - x) \right. \right.$$

$$\left. \left. \times (I(Z_i = z_1) + I(Z_i = z_2)) \right\} \right].$$

Using polynomial approximation in the neighborhood of $x$, we obtain the following approximation to $L_n$:

$$\tilde{L}_n(\rho, \beta_1^*, \beta_2^*)$$

$$= \sum_{j=1}^{N} K_h(X_{(j)} - x)[I_{1(j)}(\psi_1(x) + \widetilde{X}_{(j)}^{*T}\beta_1^*) + I_{2(j)}(\psi_2(x) + \widetilde{X}_{(j)}^{*T}\beta_2^*)]$$

$$- \sum_{j=1}^{N} K_h(X_{(j)} - x)(I_{1(j)} + I_{2(j)})$$

$$\times \log \left\{ \sum_{i \in R_j} K_h(X_i - x)[I_{1i} \exp(\psi_1(x) + \widetilde{X}_i^{*T}\beta_1^*) \right.$$

(2.14)

$$\left. + I_{2i} \exp(\psi_2(x) + \widetilde{X}_i^{*T}\beta_2^*)] \right\}$$

$$= \sum_{j=1}^{N} K_h(X_{(j)} - x)[I_{1(j)}\widetilde{X}_{(j)}^{*T}\beta_1^* + I_{2(j)}(\rho + \widetilde{X}_{(j)}^{*T}\beta_2^*)]$$

$$-\sum_{j=1}^{N} K_h(X_{(j)} - x)(I_{1(j)} + I_{2(j)})$$

$$\times \log\left\{\sum_{i \in R_j} K_h(X_i - x)[I_{1i} \exp(\widetilde{X}_i^{*T}\beta_1^*) + I_{2i} \exp(\rho + \widetilde{X}_i^{*T}\beta_2^*)]\right\},$$

where $\rho = \psi_2(x) - \psi_1(x)$, $I_{ki} = I(Z_i = z_k)$ for $k = 1, 2$, $i = 1, 2, \ldots, n$ and

$$\psi_k(x) = \psi(x, z_k), \qquad \beta_k^* = (\psi_k'(x), \psi_k''(x)/2, \ldots, \psi_k^{(p)}(x)/p) \qquad \text{for } k = 1, 2,$$

$$\widetilde{X}_i^* = ((X_i - x), (X_i - x)^2, \ldots, (X_i - x)^p) \qquad \text{for } i = 1, 2, \ldots, n.$$

By following an argument similar to the one in the previous section, we also adopt a two-step strategy here. In the first step, we apply the procedure of [12] to estimate $\beta_1^*$ and $\beta_2^*$ using observations in the neighborhoods of $x$ for the two groups separately. In the second step, we estimate $\rho$ by maximizing $\widetilde{L}_n(\rho, \hat{\beta}_1^*, \hat{\beta}_2^*)$. We now present the following theorem.

THEOREM 2. *Let $\hat{\rho}$ be the maximizer of $\widetilde{L}_n(\rho, \hat{\beta}_1^*, \hat{\beta}_2^*)$. Under conditions 1–5 of Section* 2.1, *with conditions* 2–4 *modified slightly so that they hold at point $x$, we have*

(2.15)
$$\sqrt{nh}(\hat{\rho} - \rho - b_{1n}(x)) \xrightarrow{D} N(0, \tilde{\sigma}^2(x)),$$

*where*

$$b_{1n}(x) = \frac{h^{p+1}}{(p+1)!}(\psi_2^{(p+1)}(x) - \psi_1^{(p+1)}(x))\left(\int u^{p+1} K(u)\, du\right),$$

$$\tilde{\sigma}^2(x) = \frac{(\int K^2(u)\, du)}{f(x)}\left(\frac{\sigma_1^2(x)}{p_{1x}} + \frac{\sigma_2^2(x)}{p_{2x}}\right),$$

*with $p_{kx} = P(Z = z_k | X = x)$ and, for $k = 1, 2$,*

$$\sigma_k^2 = \left[\exp(\psi_k(x))\int P(Y \geq v | X = x, Z = z_k)\lambda_0(v)\, dv\right]^{-1}$$

$$= [E\{\delta | X = x, Z = z_k\}]^{-1}.$$

Similar to the estimation of the bias term for Theorem 1, the only unknown terms $\psi_k^{(p+1)}(x)$ for $k = 1, 2$ are already obtained during the first step when we estimate the derivatives using [12] if we choose $p_1 > p$. To estimate the variance term, we note that $n^{-1}\sum_{i=1}^{n} I_{ki}\delta_i K_h(X_i - x)$ is an unbiased estimator of $p_{kx}f(x)/\sigma_k^2(x)$ for $k = 1, 2$. Natural candidates for the estimation of bias and variance are, respectively,

$$\hat{b}_{1n}(x) = \frac{h^{p+1}}{(p+1)!}(\hat{\psi}_2^{(p+1)}(x) - \hat{\psi}_1^{(p+1)}(x))\left(\int u^{p+1} K(u)\, du\right),$$

$$(2.16) \quad \hat{\sigma}^2(x) = \left(\int K^2(u)\,du\right)\left\{\frac{1}{n^{-1}\sum_{i=1}^n I_{1i}\delta_i K_h(X_i - x)} + \frac{1}{n^{-1}\sum_{i=1}^n I_{2i}\delta_i K_h(X_i - x)}\right\}.$$

As a consequence of (2.15), the theoretical optimal bandwidth minimizes the asymptotic weighted mean integrated squared error,

$$\int \left[\left\{\frac{h^{p+1}}{(p+1)!}(\psi_2^{(p+1)}(x) - \psi_1^{(p+1)}(x))\left(\int u^{p+1}K(u)\,du\right)\right\}^2 + \frac{(\int K^2(u)\,du)}{nhf(x)}\left(\frac{\sigma_1^2(x)}{p_{1x}} + \frac{\sigma_2^2(x)}{p_{2x}}\right)\right]w(x)\,dx,$$

with some weight function $w \geq 0$. We find that the asymptotically optimal constant bandwidth is given by

$$h_{\mathrm{opt}} = C_{0,p}(K) \\ \times \left[\frac{\int (1/f(x))(\sigma_1^2(x)/p_{1x} + \sigma_2^2(x)/p_{2x})w(x)\,dx}{\int (\psi_2^{(p+1)}(x) - \psi_1^{(p+1)}(x))^2 w(x)\,dx}\right]^{1/(2p+3)} n^{-1/(2p+3)},$$

with $C_{0,p}(K)$ being the same as in Section 2.1.

REMARK 5. Although the risk functions for each group can only be estimated up to a constant, the difference between two groups can be identified. Based on this observation, if we are interested in estimating risk functions for $k$ groups, we need to impose only one condition, such as $\psi_1(0) = 0$, for identifiability. On the other hand, Fan, Gijbels and King [12] needed to estimate the risk functions for each group separately. Therefore, $k$ conditions, $\psi_l(0) = 0$ for $l = 1, 2, \ldots, k$, should be imposed for identification. Sometimes, this can be inappropriate, as in the case of analyzing PBC data, to be discussed in the next section.

**3. Numerical studies.** Extensive numerical studies were conducted to evaluate the new procedures and we found that the finite sample performance of our procedure is either comparable or better than that of [12]. The Epanechnikov kernel is employed in all of the simulation studies, as well as for the analysis of the real data set.

3.1. *Simulation studies on estimating relative risk.* We compare the two procedures for the following three designs with different risk functions or distributions of the covariate variable $X$:

- Design 1: $X \sim \text{Uniform}(-1, 1)$, $\psi(x) = x^3$.

TABLE 1
*Mean integrated squared errors*

| $h_1^0$ | censoring | Design 1 [12] | Design 1 new | Design 2 [12] | Design 2 new | Design 3 [12] | Design 3 new |
|---|---|---|---|---|---|---|---|
| 0.15 | 0% | 0.181 | 0.143 | 0.324 | 0.195 | 0.262 | 0.153 |
|      | 30% | 0.271 | 0.212 | 0.384 | 0.259 | 0.353 | 0.205 |
| 0.25 | 0% | 0.091 | 0.084 | 0.165 | 0.140 | 0.141 | 0.088 |
|      | 30% | 0.129 | 0.119 | 0.213 | 0.176 | 0.197 | 0.118 |
| 0.35 | 0% | 0.062 | 0.062 | 0.152 | 0.148 | 0.104 | 0.068 |
|      | 30% | 0.084 | 0.085 | 0.182 | 0.171 | 0.155 | 0.093 |

- Design 2: $X \sim \text{Uniform}(-1,1)$, $\psi(x) = x^3 + \exp(-150(x+0.3)^2) + \exp(-150(x-0.3)^2)$.
- Design 3: $\psi(x) = x^3$. Half of the $X$ are from $N(-0.6, 0.3^2)$ truncated at $-1$ and $0$, the other half from $N(0.6, 0.3^2)$ truncated at $0$ and $1$. Note the sparsity of data in the neighborhood of $0$.

The survival time is set to be $\exp(-\psi(X) + \varepsilon)$, where $\varepsilon$ is from the standard extreme-value distribution. This is justified by the well-known result on the equivalence of the proportional hazards model (1.1) to the transformation model: $\log \Lambda_0(T) = -\psi(X) + \varepsilon$, where $\varepsilon$ is a standard extreme-value random variable. In addition, the censoring variable is assumed to be uniform on $(0, c)$, where $c$ is chosen for a prespecified censoring proportion (viz., 0% and 30%). For each set of $c$ and $\psi$, we simulate 500 realizations of $\{(X_i, Y_i, \delta_i), i = 1, 2, \ldots, 300\}$. Let $h_1^0$, $h_2^0$ and $h_3^0$ denote the bandwidths used for [12] for designs 1, 2 and 3, respectively; we tried three different values for $h_1^0$, namely, $0.15, 0.25$ and $0.35$, for which their approach yields reasonable estimates. Some adjustments were made in choosing $h_2^0$ and $h_3^0$, due to some unique features involving $\psi(\cdot)$ or the distribution of $X$. We set $h_2^0 = h_1^0$ when $|x| > 0.5$ and $h_2^0 = 0.8 h_1^0$ when $|x| \leq 0.5$. The reason for a smaller bandwidth for $|x| \leq 0.5$ is similar to the idea of variable bandwidth [10]. For design 3, we set $h_3^0 = h_1^0$ when $|x| > 0.2$ and $h_3^0 = 2 h_1^0$ when $|x| \leq 0.2$. The doubling of the bandwidth in the neighborhood of zero ensures enough data in that neighborhood. The bandwidths for our procedure are always set to be $h = 0.8 h^0$ for $h^0 = h_1^0$, $h_2^0$ and $h_3^0$ for the three designs. The mean integrated squared error (MISE) of the function estimations of $\psi(\cdot) - \psi(0)$ for designs 1 and 2 and of $\psi(\cdot) - \psi(-0.6)$ for design 3 are reported in Table 1. From the table, we find that our procedure is slightly better than [12] for design 1 and has better performance for designs 2 and 3.

We consider now the biases of the estimates. The function estimates are shown in Figures 1–3 (FGK are function estimates by Fan, Gijbels and King [12]) for the case of $h_1^0 = 0.25$ and 30% censoring; similar results are observed for $h_1^0 = 0.15$ and 0.35 and are thus omitted here for brevity. Figures 1 and 2 and other unreported results suggest that the biases of the two procedures for designs 1 and 2 are comparable. For design 2, both procedures have some biases in the neighborhoods of $\pm 0.3$ due to the peaks of $\psi(\cdot)$ at these two points. Figure 3 suggests obvious bias in [12] when $x > -0.2$ and also some bias in our procedure in the neighborhood of 0. Note that for design 3, poor performance in the neighborhood of 0 is to be expected, due to the sparsity of data.

The mean squared errors at various points for 30% censoring are reported in Table 2 to facilitate a more detailed comparison between the two approaches. We now take a close look at Table 2 to understand better the advantages of our procedure for designs 2 and 3. We find that the two procedures are largely comparable in terms of mean squared errors at all the points selected for design 1. For design 2, we observe a similar pattern for $|x| \leq 0.4$. For $|x| > 0.4$, our procedure performs better than [12]. For design 3, the performance of the two procedures is comparable for $x \leq 0.2$, and our procedure outperforms [12] when $x > 0.2$.

The better performance of our procedure for designs 2 and 3 in certain regions is largely due to the way the two estimators are constructed. In the case of design 2, due to the difficulty in estimating $\psi'(\cdot)$ in the neighborhoods of the two peaks at $x = \pm 0.3$, both procedures are not expected to estimate the function well in the neighborhoods of these two points. But, for [12], the estimates of $\psi(x) - \psi(0)$ based on $\int_0^x \psi'(u)\,du$ will be adversely affected if 0.3 or $-0.3$ lies between 0 and $x$. On the other hand, our procedure performs well, provided $\psi'(\cdot)$ can be estimated well in the neighborhoods of 0 and $x$. Similar arguments also apply to design 3, for which $\psi'(\cdot)$ is not expected to be estimated well in the neighborhoods of 0 due to the sparsity of data, whereas $\psi'(x)$ can be better estimated when $x$ is close to 0.6 or $-0.6$. This explains the obvious advantage of our procedure in the estimation of $\psi(x) - \psi(-0.6)$

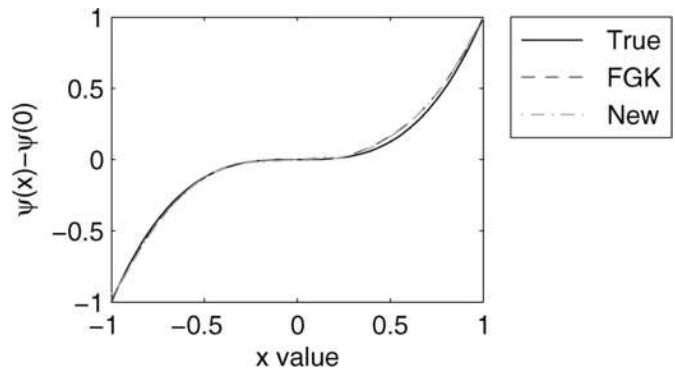

Fig. 1. *Design* 1.

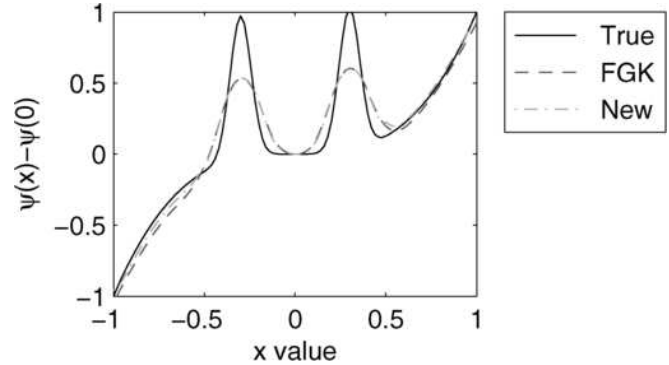

Fig. 2. *Design* 2.

in the neighborhood of $x = 0.6$ and a relatively poor performance of [12] for all $x > 0$.

3.2. *Application of estimation of the treatment effect.* We apply our procedure to the estimation of the treatment effect in an analysis of the Primary Biliary Cirrhosis (PBC) data set. Our procedure offers a natural approach to this particular problem. Basically, we estimate the treatment effect using data with bilirubin values in a neighborhood of each point of estimation. A detailed description of the PBC data can be found in Chapter 4 of [13]. A total of 312 patients participated in the randomized trial. Of the randomized patients, 187 cases (60%) were censored.

Fan and Gijbels [11] investigated the effect of treatment differences by dividing the data into two groups according to the treatment code. For each treatment group, model (2.4) was fitted using log(Bilirubin) as a covariate. The resulting curves are reproduced in Figure 4 for the sake of comparison. It is worth pointing out that, with [12], treatment differences can be estimated only up to a constant. Figure 4 implicitly assumes that the risk functions of the two treatment groups are equal at the left endpoint of the support of the covariate. There is no justification for this assumption since, while the risk functions themselves are not identifiable, the difference can be estimated following our approach. Based on Figure 4, Fan and Gijbels [11] suggested that the treatment effect is present.

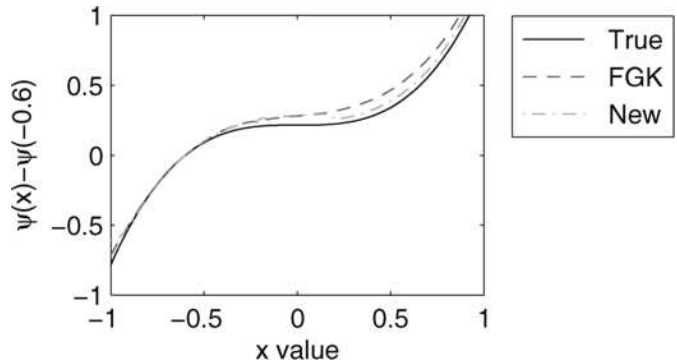

Fig. 3. *Design* 3.

TABLE 2
*Mean squared errors with 30% censoring*

|  | $x$ | $-1$ | $-0.8$ | $-0.6$ | $-0.4$ | $-0.2$ | $0$ | $0.2$ | $0.4$ | $0.6$ | $0.8$ | $1$ |
|---|---|---|---|---|---|---|---|---|---|---|---|---|
| Design 1 | | | | | | | | | | | | |
| $h_1^0 = 0.15$ | [12] | 0.745 | 0.410 | 0.363 | 0.350 | 0.338 | 0.000 | 0.341 | 0.349 | 0.370 | 0.359 | 0.602 |
| | new | 0.652 | 0.363 | 0.327 | 0.338 | 0.292 | 0.000 | 0.338 | 0.334 | 0.348 | 0.328 | 0.473 |
| $h_1^0 = 0.25$ | [12] | 0.489 | 0.272 | 0.255 | 0.252 | 0.198 | 0.000 | 0.207 | 0.264 | 0.253 | 0.255 | 0.425 |
| | new | 0.471 | 0.262 | 0.251 | 0.255 | 0.197 | 0.000 | 0.215 | 0.261 | 0.248 | 0.245 | 0.380 |
| $h_1^0 = 0.35$ | [12] | 0.406 | 0.241 | 0.210 | 0.186 | 0.121 | 0.000 | 0.125 | 0.192 | 0.200 | 0.216 | 0.356 |
| | new | 0.424 | 0.244 | 0.211 | 0.192 | 0.106 | 0.000 | 0.132 | 0.196 | 0.201 | 0.214 | 0.348 |
| Design 2 | | | | | | | | | | | | |
| $h_1^0 = 0.15$ | [12] | 0.837 | 0.484 | 0.448 | 0.451 | 0.418 | 0.000 | 0.403 | 0.429 | 0.430 | 0.421 | 0.650 |
| | new | 0.681 | 0.366 | 0.328 | 0.382 | 0.421 | 0.000 | 0.386 | 0.379 | 0.333 | 0.312 | 0.487 |
| $h_1^0 = 0.25$ | [12] | 0.529 | 0.332 | 0.321 | 0.305 | 0.279 | 0.000 | 0.277 | 0.311 | 0.307 | 0.312 | 0.438 |
| | new | 0.496 | 0.269 | 0.263 | 0.288 | 0.284 | 0.000 | 0.271 | 0.302 | 0.257 | 0.265 | 0.388 |
| $h_1^0 = 0.35$ | [12] | 0.424 | 0.273 | 0.243 | 0.232 | 0.176 | 0.000 | 0.188 | 0.249 | 0.242 | 0.243 | 0.387 |
| | new | 0.401 | 0.258 | 0.224 | 0.237 | 0.181 | 0.000 | 0.192 | 0.257 | 0.216 | 0.222 | 0.365 |
| Design 3 | | | | | | | | | | | | |
| $h_1^0 = 0.15$ | [12] | 0.758 | 0.302 | 0.000 | 0.298 | 0.388 | 0.401 | 0.366 | 0.502 | 0.467 | 0.488 | 0.707 |
| | new | 0.620 | 0.308 | 0.000 | 0.296 | 0.267 | 0.300 | 0.390 | 0.292 | 0.265 | 0.300 | 0.511 |
| $h_1^0 = 0.25$ | [12] | 0.490 | 0.194 | 0.000 | 0.178 | 0.295 | 0.311 | 0.288 | 0.383 | 0.381 | 0.372 | 0.527 |
| | new | 0.489 | 0.198 | 0.000 | 0.180 | 0.170 | 0.201 | 0.295 | 0.232 | 0.209 | 0.218 | 0.371 |
| $h_1^0 = 0.35$ | [12] | 0.424 | 0.149 | 0.000 | 0.137 | 0.249 | 0.268 | 0.270 | 0.330 | 0.371 | 0.382 | 0.461 |
| | new | 0.444 | 0.159 | 0.000 | 0.140 | 0.117 | 0.144 | 0.242 | 0.206 | 0.189 | 0.206 | 0.314 |

Following [11], we take the time (in days) between registration and death, or the time to being censored (liver transplantation or alive at study analysis) as the response and the natural logarithm of Serum Bilirubin (in mg/dl) as the continuous covariate. We use the local partial likelihood method (2.14) with $p = 1$. The derivatives $\beta_1^*$ and $\beta_2^*$ are estimated separately using the approach of [12] with $p_1 = 2$ and the bandwidth $h_1 = 1.2$, which is the same as the bandwidth used to produce Figure 5.9 of [11]. For our second-step estimator, the bandwidth is chosen to be $h = 0.8h_1$. Our 95% confidence interval for $\rho$ is constructed as

$$\hat{\rho} - \hat{b}_{1n}^A(x) \pm \Phi^{-1}(1 - \alpha/2)\hat{\sigma}(x),$$

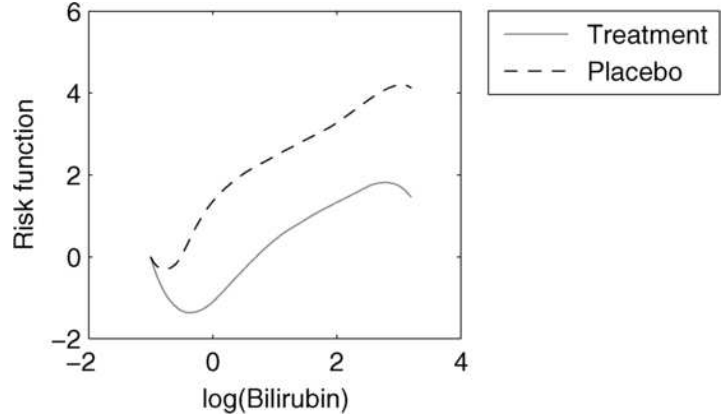

Fig. 4. *Estimates by Fan, Farmen and Gijbels.*

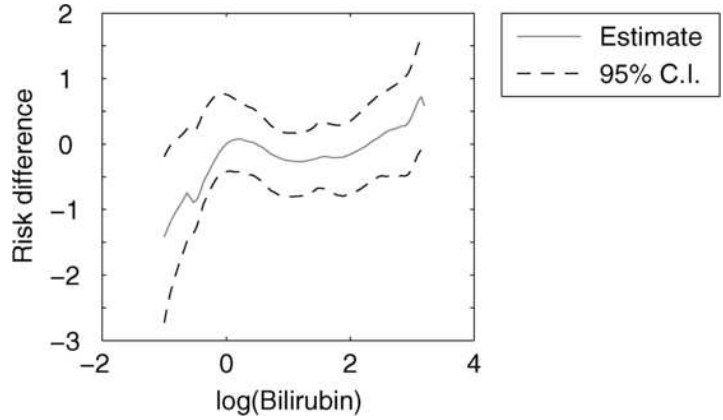

Fig. 5. *Estimates and 95% C.I.s.*

where $\hat{\tilde{b}}_{1n}(y)$ and $\hat{\sigma}(x)$ are defined in (2.16) and

$$\hat{b}_{1n}^A(x) = \int \hat{\tilde{b}}_{1n}(y) K_h(y-x)\, dx$$

is a local weighted average of estimated bias. The main reason for the average is to stabilize the bias function, which involves an estimation of a higher-order derivative curve, from abrupt change. The idea of smoothing was also adopted by Fan, Farmen and Gijbels [9]. The results are shown in Figure 5.

Figure 5 shows that, contrary to the findings of [11], the treatment effect is not present. The only sign of a treatment effect is in the range of negative covariate values. However, a close inspection of the data set reveals that the estimation for negative covariates is very unreliable since the censoring percentage is very high, at 84%, while the censoring percentage for positive covariates is 48%.

**4. Conclusion.** In this paper, we have considered direct estimation of the relative risk function in the proportional hazards model through a new version of the local partial likelihood. Our procedure was extended to the case where discrete covariates such as treatment group indicators are also present. We found that, in estimating the relative risk function, our procedure is either comparable to or outperforms the estimator proposed by Fan, Gijbels and King [12]. We applied our procedure to estimating the treatment effect in the PBC data and found that, consistent with findings from parametric analysis, a treatment effect is not present, contrary to the findings by Fan and Gijbels [11].

## APPENDIX: PROOFS

PROOF OF THEOREM 1. Set

$$u^* = \{u, u^2, \ldots, u^p\}^T, \qquad \nu_1 = \int u^* K(u)\, du,$$

(A.0) $\quad H_1^* = \mathrm{diag}(h_1, h_1^2, \ldots, h_1^p),$

$$a_x(v) = e^{\psi(x)} f(x) S(v|x), \qquad \gamma(v) = \frac{a_{x_2}(v)}{a_{x_1}(v) + a_{x_2}(v)},$$

$$\omega_x(t) = \int_0^t a_x(v)\lambda_0(v)\, dv, \qquad \kappa(t) = \int_0^t \frac{a_{x_1}(v)a_{x_2}(v)}{a_{x_1}(v) + a_{x_2}(v)} \lambda_0(v)\, dv.$$

We first prove consistency. Define

(A.1) $\qquad N_i(t) = I\{Y_i \leq t, \delta_i = 1\} \quad \text{and} \quad Y_i(t) = I\{Y_i \geq t\}.$

Let the filtration $\mathcal{F}_{nt}$ be the statistical information accruing during the time $[0, t]$, namely

$$\mathcal{F}_{nt} = \sigma\{X_i, N_i(v), Y_i(v+), i = 1, \ldots, n, 0 \leq v \leq t\}.$$

Then, under the independent censoring scheme,

(A.2) $\qquad M_i(t) = N_i(t) - \int_0^t Y_i(v) \exp(\psi(X_i))\lambda_0(v)\, dv$

is an $\mathcal{F}_{nt}$-martingale.

By (2.6)–(2.8), $\hat{\theta}_0 = (\hat{\alpha} - \alpha)$ maximizes $l_n(\theta_0, \hat{\theta}_1, \hat{\theta}_2)$, where

$$\hat{\theta}_1 = H_1^*(\hat{\beta}_{x_1}^* - \beta_{x_1}^*), \qquad \hat{\theta}_2 = H_1^*(\hat{\beta}_{x_2}^* - \beta_{x_2}^*)$$

and

$$l_n(\theta_0, \theta_1, \theta_2) = \int_0^\infty n^{-1} \sum_{i=1}^n \{K_h(X_i - x_1)(\widetilde{X}_{1i}^{*T}\beta_{x_1}^* + U_{1i}^{*T}\theta_1)$$
$$+ K_h(X_i - x_2)(\widetilde{X}_{2i}^{*T}\beta_{x_2}^* + U_{2i}^{*T}\theta_2 + \theta_0)\} \, dN_i(v)$$
$$- \int_0^\infty \log\{nS_{n,0}(\theta_0, \theta_1, \theta_2, v)\}n^{-1}$$
$$\times \sum_{i=1}^n \{K_h(X_i - x_1) + K_h(X_i - x_2)\} \, dN_i(v).$$

Here,

$$U_{1i}^* = H_1^{*-1}\widetilde{X}_{1i}^*, \qquad U_{2i}^* = H_1^{*-1}\widetilde{X}_{2i}^*,$$

$$S_{n,0}(\theta_0, \theta_1, \theta_2, v) = n^{-1} \sum_{i=1}^n Y_i(v)[\exp(\widetilde{X}_{1i}^{*T}\beta_{x_1}^* + U_{1i}^{*T}\theta_1)K_h(X_i - x_1)$$
$$+ \exp(\alpha + \theta_0)$$
$$\times \exp(\widetilde{X}_{2i}^{*T}\beta_{x_2}^* + U_{2i}^{*T}\theta_2)K_h(X_i - x_2)].$$

With a slight abuse of notation, let $\hat{\theta}_0$ maximize $l_n(\theta_0, \hat{\theta}_1, \hat{\theta}_2, \tau)$, where

$$l_n(\theta_0, \theta_1, \theta_2, \tau) = \int_0^\tau n^{-1} \sum_{i=1}^n \{K_h(X_i - x_1)(\widetilde{X}_{1i}^{*T}\beta_{x_1}^* + U_{1i}^{*T}\theta_1)$$
$$+ K_h(X_i - x_2)(\widetilde{X}_{2i}^{*T}\beta_{x_2}^* + U_{2i}^{*T}\theta_2 + \theta_0)\} \, dN_i(v)$$
(A.3)
$$- \int_0^\tau \log\{nS_{n,0}(\theta_0, \theta_1, \theta_2, v)\}n^{-1}$$
$$\times \sum_{i=1}^n \{K_h(X_i - x_1) + K_h(X_i - x_2)\} \, dN_i(v).$$

Our case corresponds to that of $\tau = \infty$.

Since $\hat{\theta}_1 \xrightarrow{p} 0, \hat{\theta}_2 \xrightarrow{p} 0$, similar to (6.26) in [12], we can show that, for any $\theta_0$,

$$l_n(\theta_0, \hat{\theta}_1, \hat{\theta}_2, \tau) = l_n(\theta_0, 0, 0, \tau) + o_p(1).$$

Let

$$S_{kn}(v) = n^{-1} \sum_{j=1}^n Y_i(v) \exp(\psi(X_i))K_h(X_i - x_k) \qquad \text{for } k = 1, 2.$$

Then

$$l_n(\theta_0, 0, 0, \tau) - l_n(0, 0, 0, \tau) = A_n(\theta_0, 0, 0, \tau) + X_n(\theta_0, 0, 0, \tau),$$

where

$$A_n(\theta_0, 0, 0, \tau) = \int_0^\tau S_{2n}(v)\theta_0\lambda_0(u)\, du$$
$$- \int_0^\tau \log\left\{\frac{S_{n,0}(\theta_0, 0, 0, v)}{S_{n,0}(0, 0, 0, v)}\right\}(S_{1n}(v) + S_{2n}(v))\lambda_0(v)\, dv,$$

$$X_n(\theta_0, 0, 0, \tau) = \int_0^\tau n^{-1}\sum_{i=1}^n K_h(X_i - x_2)\theta_0\, dM_i(v)$$
$$- \int_0^\tau \log\left\{\frac{S_{n,0}(\theta_0, 0, 0, v)}{S_{n,0}(0, 0, 0, v)}\right\} n^{-1}$$
$$\times \sum_{i=1}^n \{K_h(X_i - x_1) + K_h(X_i - x_2)\}\, dM_i(v).$$

The process $X_n(\theta_0, 0, 0, \cdot)$ is a locally integrable martingale with a predictable variation process

$$B_n(t) = \langle X_n(\theta_0, 0, 0, t), X_n(\theta_0, 0, 0, t)\rangle$$
$$= \sum_{i=1}^n \int_0^t n^{-2}\left[\left(K_h(X_i - x_2)\theta_0\right.\right.$$
$$\left.\left. - (K_h(X_i - x_1) + K_h(X_i - x_2))\log\frac{S_{n,1}(\theta_0, 0, 0, v)}{S_{n,0}(0, 0, 0, v)}\right)\right]^2$$
$$\times Y_i(v)\exp(\psi(X_i))\lambda_0(v)\, dv.$$

By Lemma 1 of [12], it can be shown that

$$EX_n^2(\theta_0, 0, 0, t) = EB_n(t) = O(n^{-1}h^{-1}) \qquad \text{for } 0 \le t \le \tau,$$

and

$$S_{n,0}(\theta_0, 0, 0, v) = S(v|x_1)f(x_1) + \exp(\alpha + \theta_0)S(v|x_2)f(x_2) + o_p(1)$$
$$= \exp(-\psi(x_1))[a_{x_1}(v) + \exp(\theta_0)a_{x_2}(v)] + o_p(1),$$
$$S_{1n}(v) = a_{x_1}(v) + o_p(1), \qquad S_{2n}(v) = a_{x_2}(v) + o_p(1).$$

Therefore,

$$A_n(\theta_0, 0, 0, \tau) = A(\theta_0, 0, 0, \tau) + o_p(1),$$

where

$$A(\theta_0, 0, 0, \tau) = \left(\int_0^\tau a_{x_2}(v)\lambda_0(v)\, dv\right)\theta_0$$
$$- \int_0^\tau \log\left\{\frac{a_{x_1}(v) + \exp(\theta_0)a_{x_2}(v)}{a_{x_1}(v) + a_{x_2}(v)}\right\}(a_{x_1}(v) + a_{x_2}(v))\lambda_0(v)\, dv.$$

Consequently,

$$l_n(\theta_0, \hat{\theta}_1, \hat{\theta}_2, \tau) = A(\theta_0, 0, 0, \tau) + o_p(1).$$

It can be easily shown that $A(\theta_0, 0, 0, \tau)$ is a strictly concave function with a unique maximum at $\theta_0 = 0$. As $l_n(\theta_0, \hat{\theta}_1, \hat{\theta}_2, \tau)$ is a concave function of $\theta_0$, by the convexity lemma [1], $\hat{\theta}_0 \xrightarrow{p} 0$.

We now consider asymptotic normality. Note that

$$
\begin{aligned}
0 &= \frac{\partial l_n(\hat{\theta}_0, \hat{\theta}_1, \hat{\theta}_2, \tau)}{\partial \theta_0} \\
&= \frac{\partial l_n(0, 0, 0, \tau)}{\partial \theta_0} + \frac{\partial^2 l_n(\bar{\theta}_0, \bar{\theta}_1, \bar{\theta}_2, \tau)}{\partial \theta_0^2} \hat{\theta}_0 \\
&\quad + \frac{\partial^2 l_n(\bar{\theta}_0, \bar{\theta}_1, \bar{\theta}_2, \tau)}{\partial \theta_0 \, \partial \theta_1} \hat{\theta}_1 + \frac{\partial^2 l_n(\bar{\theta}_0, \bar{\theta}_1, \bar{\theta}_2, \tau)}{\partial \theta_0 \, \partial \theta_2} \hat{\theta}_2,
\end{aligned}
\quad (\text{A.4})
$$

where $\bar{\theta}_k$ lies on the line segment from 0 to $\hat{\theta}_k$ for $k = 0, 1, 2$. By making use of Lemma 1 from [12], it is straightforward (although tedious) to show that

$$
\begin{aligned}
\frac{\partial^2 l_n(\bar{\theta}_0, \bar{\theta}_1, \bar{\theta}_2, \tau)}{\partial \theta_0^2} &= -\kappa(\tau) + o_p(1), \\
\frac{\partial^2 l_n(\bar{\theta}_0, \bar{\theta}_1, \bar{\theta}_2, \tau)}{\partial \theta_0 \, \partial \theta_1^T} &= \kappa(\tau)\nu_1^T + o_p(1), \\
\frac{\partial^2 l_n(\bar{\theta}_0, \bar{\theta}_1, \bar{\theta}_2, \tau)}{\partial \theta_0 \, \partial \theta_2^T} &= -\kappa(\tau)\nu_1^T + o_p(1),
\end{aligned}
\quad (\text{A.5})
$$

where $\kappa(\tau)$ and $\nu_1$ are defined in (A.0). Furthermore, we prove at the end of the Appendix that

$$\frac{\partial l_n(0, 0, 0, \tau)}{\partial \theta_0} = U_{0n}(\tau) + b_{0n}(\tau) + o_p(h^{p+1}), \quad (\text{A.6})$$

where

$$
\begin{aligned}
U_{0n}(\tau) = n^{-1} \sum_{i=1}^{n} \int_0^{\tau} &\bigg\{ K_h(X_i - x_2) \\
&\quad - \frac{S_{n,1}(0, 0, 0, v)}{S_{n,0}(0, 0, 0, v)} [K_h(X_i - x_1) + K_h(X_i - x_2)] \bigg\} dM_i(v),
\end{aligned}
\quad (\text{A.7})
$$

$$b_{0n}(\tau) = \frac{h^{p+1}\kappa(\tau)}{(p+1)!} \{\psi^{(p+1)}(x_2) - \psi^{(p+1)}(x_1)\} \left( \int u^{p+1} K(u) \, du \right).$$

Here,

$$S_{n,1}(\theta_0, \theta_1, \theta_2, v)$$

$$= \frac{\partial S_{n,0}(\theta_0, \theta_1, \theta_2, v)}{\partial \theta_0}$$

$$= n^{-1} \sum_{j=1}^{n} Y_i(v) \exp(\alpha^0 + \theta_0) \exp(\widetilde{X}_{2i}^{*T} \beta_{x_2}^{*0} + U_{2i}^{*T} \theta_2) K_h(X_i - x_2).$$

It follows from (A.4)–(A.6) that

(A.8)
$$\hat{\theta}_0 = (1 + o_p(1)) \left[ \frac{U_{0n}(\tau)}{\kappa(\tau)} + \frac{b_{0n}(\tau)}{\kappa(\tau)} \right.$$
$$\left. + (\nu_1 + o_p(1))^T (\hat{\theta}_1 - \hat{\theta}_2) + o_p(h^{p+1}) \right].$$

By applying the martingale property [13], we can easily prove that

(A.9) $$\sqrt{nh} U_{0n}(\tau) \xrightarrow{d} N\left(0, \kappa(\tau) \int K^2(u) \, du\right).$$

In addition, it follows from Theorem 4 of [12] that $\sqrt{nh_1} \hat{\theta}_l = O_p(1)$ for $l = 1, 2$. We conclude that

(A.10) $$\sqrt{nh} \hat{\theta}_l = \sqrt{\frac{h}{h_1}} \sqrt{nh_1} \hat{\theta}_l = o_p(1), \qquad \text{for } l = 1, 2,$$

since $h/h_1 \to 0$ by Condition 5 of Theorem 1. Finally, from (A.7)–(A.10), we have

$$\sqrt{nh}(\hat{\theta}_0 - b_n(\tau)) \xrightarrow{d} N\left(0, \frac{(\int K^2(u) \, du)}{\kappa(\tau)}\right),$$

where

$$b_n(\tau) = \frac{b_{0n}(\tau)}{\kappa(\tau)} = \frac{h^{p+1}}{(p+1)!} \{\psi^{(p+1)}(x_2) - \psi^{(p+1)}(x_1)\} \left(\int u^{p+1} K(u) \, du\right).$$

To finish the proof of Theorem 1, it suffices to prove (A.6). □

PROOF OF (A.6). By taking the derivative with respect to $\theta_0$ in (A.1), we obtain

$$\frac{\partial l_n(0,0,0,\tau)}{\partial \theta_0} = U_{0n}(\tau) + B_{0n}(\tau),$$

where

$$U_{0n}(\tau) = n^{-1} \sum_{i=1}^{n} \int_0^{\tau} \left\{ K_h(X_i - x_2) \right.$$
$$\left. - \frac{S_{n,1}(0,0,0,v)}{S_{n,0}(0,0,0,v)} [K_h(X_i - x_1) + K_h(X_i - x_2)] \right\} dM_i(v)$$

and

$$B_{0n}(\tau) = n^{-1} \sum_{i=1}^{n} \int_0^{\tau} \bigg\{ K_h(X_i - x_2)$$
$$- \frac{S_{n,1}(0,0,0,u)}{S_{n,0}(0,0,0,u)}[K_h(X_i - x_1) + K_h(X_i - x_2)]\bigg\}$$
$$\times Y_i(u)\exp(\psi(X_i))\lambda_0(u)\,du$$
$$= B_{01n}(\tau) + B_{02n}(\tau).$$

Here

$$B_{01n}(\tau) = -n^{-1} \sum_{i=1}^{n} \int_0^{\tau} K_h(X_i - x_1) \frac{S_{n,1}(0,0,0,v)}{S_{n,0}(0,0,0,v)} Y_i(v)\exp(\psi(X_i))\lambda_0(v)\,dv$$

and

$$B_{02n}(\tau) = n^{-1} \sum_{i=1}^{n} \int_0^{\tau} K_h(X_i - x_2)$$
$$\times \bigg\{ 1 - \frac{S_{n,1}(0,0,0,v)}{S_{n,0}(0,0,0,v)}\bigg\} Y_i(v)\exp(\psi(X_i))\lambda_0(u)\,dv.$$

Define

$$B_{01n}^e(\tau) = -n^{-1} \sum_{i=1}^{n} \int_0^{\tau} K_h(X_i - x_1) \frac{S_{n,1}(0,0,0,v)}{S_{n,0}(0,0,0,v)}$$
$$\times Y_i(v)\exp(\psi(x_1) + \widetilde{X}_{1i}^{*T}\beta_{x_1}^{*0})\lambda_0(v)\,dv,$$
$$B_{02n}^e(\tau) = n^{-1} \sum_{i=1}^{n} \int_0^{\tau} K_h(X_i - x_2)\bigg\{1 - \frac{S_{n,1}(0,0,0,v)}{S_{n,0}(0,0,0,v)}\bigg\}$$
$$\times Y_i(v)\exp(\psi(x_2) + \widetilde{X}_{2i}^{*T}\beta_{x_2}^{*0})\lambda_0(v)\,dv.$$

Similar to the proof of (6.24) in [12], and by the fact that

$$\frac{S_{n,1}(0,0,0,u)}{S_{n,0}(0,0,0,u)} \to \frac{\exp(\alpha^0)f(x_2)P(v|x_2)}{f(x_1)P(u|x_1) + \exp(\alpha^0)f(x_2)P(v|x_2)}$$
$$= \frac{a_{x_2}(v)}{a_{x_1}(v) + a_{x_2}(v)} = \gamma(v),$$

we have

$$B_{01n}(\tau) - B_{01n}^e(\tau) = -\frac{\psi^{(p+1)}(x_1)h^{p+1}}{(p+1)!}\varkappa(\tau)\int u^{p+1}K(u)\,du + o_p(h^{p+1}),$$

$$B_{02n}(\tau) - B_{02n}^e(\tau) = \frac{\psi^{(p+1)}(x_2)h^{p+1}}{(p+1)!}\varkappa(\tau)\int u^{p+1}K(u)\,du + o_p(h^{p+1}).$$

Furthermore,

$$B_{02n}^e(\tau) + B_{01n}^e(\tau)$$
$$= n^{-1} \sum_{i=1}^{n} \int_0^\tau Y_i(v) K_h(X_i - x_2) \exp(\psi(x_2) + \widetilde{X}_{2i}^{*T} \beta_{x_2}^{*0}) \lambda_0(v) \, dv$$
$$- \int_0^\tau \frac{S_{n,1}(0,0,0,v)}{S_{n,0}(0,0,0,v)} n^{-1}$$
$$\times \sum_{i=1}^{n} Y_i(v) \{K_h(X_i - x_1) \exp(\psi(x_1) + X_{1i}^{*} \beta_{x_1}^{*0})$$
$$+ K_h(X_i - x_2) \exp(\psi(x_2) + \widetilde{X}_{2i}^{*T} \beta_{x_2}^{*0})\} \lambda_0(v) \, dv$$
$$= \exp(\psi(x_1)) \Big\{ \int_0^\tau S_{n,1}(0,0,0,u) \lambda_0(u) \, du$$
$$- \int_0^\tau \frac{S_{n,1}(0,0,0,u)}{S_{n,0}(0,0,0,u)} S_{n,0}(0,0,0,u) \lambda_0(u) \, du \Big\} = 0.$$

We conclude that

$$B_{0n}(\tau) = (B_{01n}(\tau) - B_{01n}^e(\tau)) + (B_{02n}(\tau) - B_{02n}^e(\tau)) = b_{0n}(\tau) + o_p(h^{p+1}),$$

where $b_{0n}(\tau)$ is defined in (A.7). □

PROOF OF THEOREM 2. Define

$$\Lambda(\tau, x) = \int_0^\tau S(v|x) \lambda_0(v) \, dv,$$

(A.11)
$$\iota(x) = \frac{p_{2x} \exp(\psi_2(x))}{p_{1x} \exp(\psi_1(x)) + p_{2x} \exp(\psi_2(x))},$$

$$\tilde{\kappa}(x) = \frac{p_{1x} \exp(\psi_1(x)) p_{2x} \exp(\psi_2(x))}{p_{1x} \exp(\psi_1(x)) + p_{2x} \exp(\psi_2(x))}.$$

The proof of Theorem 2 is similar to that of Theorem 1. Therefore, we will only provide the main steps that differ from those in the proof of Theorem 1. Denote $\hat{\eta} = (\hat{\eta}_0, \hat{\eta}_1, \hat{\eta}_2)$ with $\hat{\eta}_0 = \hat{\rho} - \rho^0$, $\hat{\eta}_k = H_1^*(\hat{\beta}_k^* - \beta_k^{*0})$ for $k = 1, 2$. Then, maximizing $\widetilde{L}_n$ of Section 2.2 is equivalent to maximizing $l_n(\eta, \infty)$, where

(A.12)
$$l_n(\eta, \tau) = \int_0^\tau n^{-1} \sum_{i=1}^{n} K_h(X_i - x)$$
$$\times [I_{1i}(\widetilde{X}_i^{*T} \beta_1^* + \tilde{U}_i^{*T} \eta_1)$$
$$+ I_{2i}(\rho + \eta_0 + \widetilde{X}_i^{*T} \beta_2^* + U_i^{*T} \eta_2)] \, dN_i(v)$$

$$- \int_0^\tau n^{-1} \sum_{i=1}^n K_h(X_i - x)(I_{1i} + I_{2i}) \log\{nS_n(\eta, v)\} \, dN_i(v).$$

Here, $N_i(t)$ [and later $M_i(t)$] are defined in (A.1) and (A.2) in the proof of Theorem 1,

$$U_i^* = (H_1^*)^{-1} \widetilde{X}_i^*, \qquad \text{with } H_1^* \text{ defined in (A.0), and}$$

$$S_n(\eta, v) = n^{-1} \sum_{i=1}^n Y_i(v) K_h(X_i - x)[I_{1i} \exp(\widetilde{X}_i^{*T} \beta_1^* + \tilde{U}_i^{*T} \eta_1)$$
$$+ I_{2i} \exp(\rho + \eta_0 + \widetilde{X}_i^{*T} \beta_2^* + U_i^{*T} \eta_2)].$$

By an argument similar to that in the proof of Theorem 1, and by noting that

$$S_n(\eta_0, 0, 0, v)$$
$$\to p_{1x} \exp(\psi_1(x)) + p_{2x} \exp(\psi_2(x) + \eta_0) \exp(-\psi_1(x)) f(x) S(v|x),$$

we can prove that

$$l_n(\eta_0, \hat{\eta}_1, \hat{\eta}_2, \tau) - l_n(0, 0, 0, \tau)$$
$$= (p_{1x} \exp(\psi_1(x)) + p_{2x} \exp(\psi_2(x)) f(x) \Lambda(\tau, x))$$
$$\times \left[\eta_0 - \log\left\{\frac{p_{1x} \exp(\psi_1(x)) + p_{2x} \exp(\psi_2(x) + \eta_0)}{p_{1x} \exp(\psi_1(x)) + p_{2x} \exp(\psi_2(x))}\right\}\right] + o_p(1).$$

Obviously, the right-hand side of the above equation is a strictly concave function of $\eta_0$. Since $l_n(\eta_0, \hat{\eta}_1, \hat{\eta}_2, \tau)$ is a concave function of $\eta_0$, by the convexity lemma, $\hat{\eta} \xrightarrow{p} 0$.

Next, we prove asymptotic normality. Note that

(A.13)
$$0 = \left.\frac{\partial l_n(\eta, \tau)}{\partial \eta_0}\right|_{\eta = \hat{\eta}}$$
$$= \frac{\partial l_n(0, \tau)}{\partial \eta_0} + \frac{\partial^2 l_n(\bar{\eta}, \tau)}{\partial \eta_0^2} \hat{\eta}_0 + \frac{\partial^2 l_n(\bar{\eta}, \tau)}{\partial \eta_0 \, \partial \eta_1^T} \hat{\eta}_1 + \frac{\partial^2 l_n(\bar{\eta}, \tau)}{\partial \eta_0 \, \partial \eta_2^T} \hat{\eta}_2,$$

where $\bar{\eta}$ lies on the line segment between 0 and $\hat{\eta}$. Using arguments similar to those in the proof of Theorem 1, we can show that for any $\bar{\eta} \xrightarrow{p} 0$,

(A.14)
$$\left.\frac{\partial^2 l_n(\eta, \tau)}{\partial \eta_0^2}\right|_{\eta = \bar{\eta}} = -\tilde{\kappa}(x) f(x) \Lambda(\tau, x) + o_p(1),$$
$$\left.\frac{\partial^2 l_n(\eta, \tau)}{\partial \eta_0 \, \partial \eta_1^T}\right|_{\eta = \bar{\eta}} = \tilde{\kappa}(x) f(x) \Lambda(\tau, x) \nu_1^T + o_p(1),$$
$$\left.\frac{\partial^2 l_n(\eta, \tau)}{\partial \eta_0 \, \partial \eta_2^T}\right|_{\eta = \bar{\eta}} = -\tilde{\kappa}(x) f(x) \Lambda(\tau, x) \nu_1^T + o_p(1),$$

with $\tilde{\kappa}(x)$ and $\Lambda(\tau, x)$ defined in (A.11). $\frac{\partial l_n(0,\tau)}{\partial \eta_0}$ can be expressed as

$$\frac{\partial l_n(0,\tau)}{\partial \eta_0} = \tilde{D}(\tau) + \tilde{b}_{1n}(\tau), \tag{A.15}$$

where

$$\tilde{D}_n(\tau) = \int_0^\tau n^{-1} \sum_{i=1}^n K_h(X_i - x)\{I_{2i} - (I_{1i} + I_{2i})q_n(v)\} \, dM_i(v),$$

$$\tilde{b}_{1n}(\tau) = \int_0^\tau n^{-1} \sum_{i=1}^n K_h(X_i - x)\{I_{2i} - (I_{1i} + I_{2i})q_n(v)\}$$
$$\times Y_i(v) \exp(\psi(X_i, Z_i))\lambda_0(v) \, dv,$$

with

$$q_n(v) = \frac{n^{-1} \sum_{i=1}^n Y_i(v) K_h(X_i - x) I_{2i} \exp(\rho^0 + \eta_0 + \widetilde{X}_i^{*T} \beta_2^{0*})}{S_n(0, v)}.$$

By Taylor expanding at $(x, z_1)$ and $(x, z_2)$,

$$I_{1i}(\exp \psi(X_i, Z_i) - \exp(\psi_1(x) + \widetilde{X}_i^{*T} \beta_1^*))$$
$$= I_{1i} \exp(\psi_1(x)) \frac{\psi_1^{(p+1)}(x)}{(p+1)!}(X_i - x)^{p+1} + o_p(h^{p+1}),$$

$$I_{2i}(\exp \psi(X_i, Z_i) - \exp(\psi_2(x) + \widetilde{X}_i^{*T} \beta_2^*))$$
$$= I_{2i} \exp(\psi_2(x)) \frac{\psi_2^{(p+1)}(x)}{(p+1)!}(X_i - x)^{p+1} + o_p(h^{p+1}),$$

we obtain for $\tilde{b}_{1n}(\tau)$

$$\tilde{b}_{1n}(\tau) = \frac{h^{p+1}}{(p+1)!}(\psi_2^{(p+1)}(x) - \psi_1^{(p+1)}(x))$$
$$\times \tilde{\kappa}(x) f(x) \Lambda(\tau, x) \left( \int u^{p+1} K(u) \, du \right) + o_p(h^{p+1}). \tag{A.16}$$

Furthermore, since

$$\langle \sqrt{nh} \tilde{D}_n(\tau), \sqrt{nh} \tilde{D}_n(\tau) \rangle$$
$$= h \int_0^\tau \sum_{i=1}^n K_h^2(X_i - x)\{-I_{1i} q_n(v) + I_{2i}(1 - q_n(v))\}^2$$
$$\times Y_i(v) \exp(\psi(X_i, Z_i))\lambda_0(v) \, dv$$
$$= h \int_0^\tau \sum_{i=1}^n K_h^2(X_i - x)\{I_{1i} q_n^2(v) + I_{2i}(1 - q_n(v))^2\}$$

$$\times Y_i(v) \exp(\psi(X_i, Z_i)) \lambda_0(v) \, dv$$

$$\xrightarrow{p} (p_{1x} \exp(\psi_1(x)) \iota^2(x) + p_{2x} \exp(\psi_2(x))(1 - \iota(x))^2)$$

$$\times f(x) \Lambda(\tau, x) \left( \int K^2(u) \, du \right)$$

$$= \tilde{\kappa}(x) f(x) \Lambda(\tau, x) \left( \int K^2(u) \, du \right) \triangleq \sigma^2(\tau, x),$$

a straightforward application of the martingale central limit theorem results in

(A.17) $$\sqrt{nh} D_n(\tau) \xrightarrow{d} N(0, \sigma^2(\tau, x)).$$

From (A.13)–(A.17) and under the condition $h/h_1 \to 0$,

$$\sqrt{nh}(\hat{\rho} - \rho - b_n(\tau)) \xrightarrow{d} N(0, \sigma(\tau)^2),$$

where

$$b_n(\tau) = \frac{h^{p+1}}{(p+1)!} (\psi_2^{(p+1)}(x) - \psi_1^{(p+1)}(x)) \left( \int u^{p+1} K(u) \, du \right),$$

$$\sigma(\tau)^2 = \frac{(\int K^2(u) \, du)}{\tilde{\kappa}(x) f(x) \Lambda(\tau, x)}. \qquad \square$$

**Acknowledgment.** The authors thank the referees for their comments that led to improvement of the paper.

Hong Kong University Science
and Technology
Clearwater Bay
Kowloon
Hong Kong
E-mail: lzzhou@ust.hk